\documentclass[12pt]{article}
\usepackage{amssymb}
\usepackage{amsmath}
\usepackage{geometry}
\usepackage{setspace}

\input{tcilatex}

\begin{document}

\title{Uniformly Balanced Repeated Measurements Designs in the Presence of
Subject Dropout}
\author{Dibyen Majumdar, Angela M. Dean and Susan M. Lewis \and University
of Illinois at Chicago, The Ohio State University \and and University of
Southampton \and \textit{University of Southampton Technical Report 388,
Second Revision}}
\maketitle

\begin{abstract}
\setstretch{1}Low, Lewis and Prescott (1999) showed that a crossover design
based on a Williams Latin square of order $4$ can suffer substantial loss of
efficiency if some observations in the final period are unavailable. Indeed,
if all observations are missing, the design becomes disconnected. We derive
the information matrix for the direct effects of a Uniformly Balanced
Repeated Measurements Design (UBRMD) in $t$ periods when subjects may drop
out before the end of the study and examine the maximum loss of information.
The special case of loss of observations in the final period only is
examined in detail. In particular we show that a UBRMD in $t\geq 5$ periods
remains connected when some or all observations in the final period are
unavailable.

\medskip


\noindent \textbf{Running Title}: Crossover Designs with Subject
Dropout\medskip

\noindent \textbf{Addresses:}{\small \ \noindent Dibyen Majumdar
(corresponding author), Department of Math., Stat., and Comp. Sci.,
University of Illinois at Chicago, 851 S. Morgan, Chicago, Il 60607-7045,
US. Phone: 312-996-4833, Fax: 312-996-1491, E-mail: dibyen@uic.edu.}

{\small \noindent Angela M. Dean, Department of Statistics, 404 Cockins
Hall, 1958 Neil Avenue, Columbus, OH 43210-1247, US. E-mail:
amd@stat.ohio-state.edu. }

{\small \noindent Susan M. Lewis, School of Mathematics, University of
Southampton, Southampton S0171BJ, UK. E-mail: S.M.Lewis@maths.soton.ac.uk.}%
\setstretch{1}
\end{abstract}

\section{Introduction}

Cross-over experiments are widely used for comparing the responses to
various different stimuli or treatments in areas ranging from psychology and
human factor engineering to medical and agricultural applications; see, for
example, the books by Jones and Kenward (2003), Ratkowsky, Evans and
Alldredge (1992) and Senn (2002). Such experiments extend over a sequence of
time periods. Each subject receives one treatment per period and an
observation is made at the end of the period. The influence of a treatment
on the subject's response may extend (or carry over) into the period
following that in which it is administered. This is known as a \emph{%
first-order carry-over effect} or \emph{first-order residual effect}. In a
simple statistical model for crossover studies, the response for a given
subject in a given period is regarded as a sum of the effects of the
subject, the period, the treatment given in this period (the \emph{direct
effect} of the treatment), the carry-over effect from the treatment given in
the preceding period, and a random error.

There is an extensive literature that assures us that a\ carefully designed
crossover study can produce a wealth of information and that the parameters
of interest can be estimated with high precision; see, for instance, Stufken
(1996). This is based on the implicit, but critical, assumption that the
experiment yields all planned observations. Yet, in many studies such as
clinical trials, there is a substantial probability that some subjects will
drop out of the study prior to completion of their treatment sequence. Low,
Lewis and Prescott (1999) observed that a dropout rate of between 5\% and
10\% is not uncommon and, in some areas, can be as high as 25\%. They give
an example of a design in four periods, based on a Williams Latin square
(Williams (1949)), where there is substantial loss of information if some
observations are unavailable in period $4$. Indeed, if all observations in
the final period are not available, the design becomes disconnected, i.e.,
elementary contrasts are no longer all estimable.

It is important to note that a similar situation may arise in an interim
analysis. When interim results on a cross-over experiment are analyzed, the
interim design may consist of the planned design without the final several
periods.

The loss of connectedness is the most severe consequence of the
unavailability of observations. A general study of the loss of connectedness
that results from unavailability of observations was done by Ghosh (1979,
1982). For crossover designs, Low, Lewis and Prescott (1999) formulated
requirements for a planned design to be robust to dropouts in terms of the
properties of the implemented designs that might result under a
\textquotedblleft completely-at-random\textquotedblright\ dropout mechanism
(Diggle and Kenward (1994)). Godolphin (2004) also studied the problem of
loss of connectedness of various designs, including crossover designs.

An experimenter generally starts with a design, the \emph{planned design,}
that possesses desirable properties, including high efficiency or
optimality. If no subject drops out, the study yields the entire information
that was envisioned at the planning stage. Dropouts, however, lead to loss
of information. The \emph{implemented} \emph{design} is the design that
corresponds to all \textit{available} observations, and this design can be
identified only at the conclusion of the experiment.

For the case of a Williams Latin Square of order $4$ the \textit{expected}
information loss for various probability distributions of dropouts was
studied by Low, Lewis and Prescott (1999). In this article, we focus on the
\textit{maximum} \textit{information loss} that may be anticipated. For
instance, in a study with four periods where subjects are expected to remain
at least through the first three periods, dropouts, if any, would occur in
the final period only. In this case minimal information is attained when
\textit{all} subjects drop out in the final period, which gives the \emph{%
minimal design.}

In this paper we assume that the \emph{planned design} belongs to the class
of \emph{Uniform Balanced Repeated Measurement Designs} (UBRMDs). This is an
important class of designs that have been studied extensively in the
literature and are a popular choice in practice. UBRMDs have elegant \textit{%
combinatorial} \textit{balance} and, under the simple model with additive
i.i.d. errors with constant variance, possess various optimality properties;
see, for example, Hedayat and Afsarinejad (1978), Cheng and Wu (1980),
Kunert (1984), Hedayat and Yang (2003), and Hedayat and Yang (2004). (Refer
to Stufken (1996) and Hedayat and Yang (2003) for additional references). A
design is called \emph{uniform} if (a) for each subject, each treatment is
allocated to the same number of periods, and (b) for each period, each
treatment is allocated to the same number of subjects. Furthermore, a design
is called \emph{balanced\/ for carryover effects }(\emph{balanced, }in
short) if, in the order of application, each treatment is preceded by every
other treatment the same number of times and is not preceded by itself.

The goal of this research is to study the maximum loss of information and
the resulting loss of precision of the estimators that result from subject
dropout when the planned design is a UBRMD. Since the maximum loss is
attained by the \emph{minimal design}, we study properties of this design,
including its information matrix and efficiency. If the maximum loss of
information is not deemed to be large, then the experimenter may conclude
that no modification of the plan for the experiment is needed. On the other
hand, if the loss is large, the experimenter should consider alternative
strategies.

We work in the same setup as Low, Lewis and Prescott (1999); in particular
we assume a completely-at-random\ dropout mechanism. Also, we assume
throughout that a subject who leaves the study does not re-enter. In Section
2 we derive general formulae for, and study the properties of, the
information matrix of the direct effects of the minimal design when the
planned design is a UBRMD, with subject dropouts occurring in the final $m$
periods only. We examine the connectedness of the minimal design and, in
particular, show that a UBRMD based on $t\geq 5$ treatments remains
connected even when all observations in the final period are unavailable. We
also develop measures for the maximum loss of precision due to subject
dropout and the efficiency of the minimal design. In Section 3 we study the
case of one-period dropout in more detail, including the special case when
the planned design is based on Williams Latin squares. Also in this section,
we identify members of the class of UBRMDs for which the loss of information
is small.

The focus of this paper is to study certain properties of UBRMDs in the
presence of subject dropout. For the broader problem of designing a
crossover experiment in the presence of subject dropout one has to choose a
planned design from a class (not necessarily the class of UBRMDs) of highly
efficient designs for which the loss of information is small, and the
corresponding minimal design is highly efficient.

\section{Setup and general results}

Consider a planned design with $p$ periods, $s$ subjects and $t$ treatments.
The simple model for the vector of response variables obtained from the
implemented design can be written as

\begin{equation}
Y=X_{S}\beta +X_{P}\alpha +X_{D}\tau +X_{C}\rho +\epsilon ,  \label{modl}
\end{equation}

\noindent where $\epsilon $ is the vector of random error variables, $\beta $
is a vector of $s$ subject effects, $\alpha $ is a vector of $p$ time period
effects, $\tau $ is a vector of $t$ direct treatment effects, $\rho $ is a
vector of $t$ carry-over effects, and the $X$ matrices are the corresponding
design matrices. The treatments are labelled $0,1,...,t-1.$ For the purposes
of designing efficient experiments, all effects in the model are assumed to
be fixed effects.

We define the following incidence and replication matrices: $%
N_{SD}~=~X_{S}^{\prime }X_{D},$ $N_{SC}~=~X_{S}^{\prime }X_{C},$ $%
N_{PD}~=~X_{P}^{\prime }X_{D},$ $N_{PC}~=~X_{P}^{\prime }X_{C}\,,$ $%
N_{DC}=X_{D}^{\prime }X_{C},$ $r_{D}=N_{DS}\mbox{\bf{1}}_{s}=N_{DP}%
\mbox{\bf{1}}_{p},$ $r_{C}=N_{CS}\mbox{\bf{1}}_{s}=N_{CP}\mbox{\bf{1}}_{p},$
where the \textquotedblleft prime\textquotedblright\ denotes transpose and $%
\mbox{\bf{1}}_{a}$ is a vector of $a$ unit elements. Also we define $%
J_{a\times b}=\mbox{\bf{1}}_{a}\mbox{\bf{1}}_{b}^{\prime }$, $%
J_{a}=J_{a\times a}$, $N_{ji}=N_{ij}^{\prime }$ (for $i,j=S,P,D,C$), $%
r_{D}^{\delta }=diag(r_{D}),$ and $r_{C}^{\delta }=diag(r_{C}).$ Moreover, $%
I_{a}$ denotes an $a\times a$ identity matrix. We order the responses period
by period for each subject in turn, so that, $X_{P}=\mbox{\bf{1}}_{s}\otimes
I_{p}~~$and$~~X_{S}=I_{s}\otimes \mbox{\bf{1}}_{p}.$

The joint information matrix for estimating the direct and carry-over
(residual) treatment effects is given by%
\begin{eqnarray}
C &=&\left[
\begin{array}{ll}
C_{11} & C_{12} \\
C_{21} & C_{22}%
\end{array}%
\right] ,  \label{eqGENinf} \\
&&  \notag \\
\text{where }C_{11} &=&r_{D}^{\delta }+\frac{1}{ps}r_{D}r_{D}^{\prime }-%
\frac{1}{p}N_{DS}N_{SD}-\frac{1}{s}N_{DP}N_{PD}  \notag \\[2ex]
C_{22} &=&r_{C}^{\delta }+\frac{1}{ps}r_{C}r_{C}^{\prime }-\frac{1}{p}%
N_{CS}N_{SC}-\frac{1}{s}N_{CP}N_{PC}  \notag \\[2ex]
C_{12} &=&N_{DC}~+\frac{1}{ps}r_{D}r_{C}^{\prime }-~\frac{1}{p}%
N_{DS}N_{SC}~-~\frac{1}{s}N_{DP}N_{PC}.  \notag
\end{eqnarray}

\noindent The information matrices for the direct effects and the carry-over
effects, respectively, are
\begin{eqnarray}
C_{D} &=&C_{11}-C_{12}C_{22}^{-}C_{21}\,,  \label{eqCD} \\[1ex]
C_{R} &=&C_{22}-C_{21}C_{11}^{-}C_{12}\,.  \label{eqCR}
\end{eqnarray}%
In this article we focus primarily on $C_{D}$.

Throughout, we assume that the \emph{planned design} $d_{\mathrm{plan}}$ is
a UBRMD with $p=t$ time periods, $s=gt$ subjects, based on $t$ treatments
and, in the \emph{implemented design} $d_{\mathrm{imp}}$, all subjects
complete their allocated treatment sequence in the first $t-m$ periods ($%
1\leq m<t-1$). After the first $t-m$ periods, subjects may start dropping
out of the study completely at random. Since we assume that, once a subject
drops out of the study, the subject will not return, the \textit{worst case
scenario} occurs when all subjects drop out at period $t-m.$ The design $%
d_{\min }$, composed of the first $t-m$ periods of $d_{\mathrm{plan}},$ is
called the \emph{minimal design}.

For even $t,$ a Williams Latin Square gives a UBRMD, as does any
sequentially counterbalanced Latin Square (see Isaac, Dean and Ostrom
(2001), for a survey). For $t$ odd, a UBRMD cannot be constructed from one
Williams Latin square, but such a design with $2t$ subjects can be
constructed from a pair of squares. In addition, when $t$ is a composite
number, Higham (1998) has shown that there exists a UBRMD in $t$ subjects, $%
t $ periods and $t$ treatments. The union of UBRMDs (identical or distinct)
with the same value of $t$ is a UBRMD. Here are some examples.

\medskip

\noindent \textbf{Example 1} \ Three UBRMDs are shown below, where the
columns show the treatment sequences and the rows correspond to the time
periods. The design $d_{2plan}$ is a Williams Latin Square of order 4, while
designs $d_{1plan}$ and $d_{3plan}$ consist of a pair of Williams Latin
Squares for $t=3$ and $t=5$ treatments, respectively.%
\begin{equation*}
\begin{array}{ccccc}
\begin{array}{c}
d_{1plan} \\
\begin{array}{cccccc}
1 & 2 & 0 & 2 & 0 & 1 \\
0 & 1 & 2 & 0 & 1 & 2 \\
2 & 0 & 1 & 1 & 2 & 0%
\end{array}%
\end{array}
&  &
\begin{array}{c}
d_{2plan} \\
\begin{tabular}{cccc}
$0$ & $1$ & $2$ & $3$ \\
$1$ & $2$ & $3$ & $0$ \\
$3$ & $0$ & $1$ & $2$ \\
$2$ & $3$ & $0$ & $1$%
\end{tabular}%
\end{array}
&  &
\begin{array}{c}
d_{3plan} \\
\begin{tabular}{cccccccccc}
$1$ & $2$ & $3$ & $4$ & $0$ & $3$ & $4$ & $0$ & $1$ & $2$ \\
$0$ & $1$ & $2$ & $3$ & $4$ & $4$ & $0$ & $1$ & $2$ & $3$ \\
$2$ & $3$ & $4$ & $0$ & $1$ & $2$ & $3$ & $4$ & $0$ & $1$ \\
$4$ & $0$ & $1$ & $2$ & $3$ & $0$ & $1$ & $2$ & $3$ & $4$ \\
$3$ & $4$ & $0$ & $1$ & $2$ & $1$ & $2$ & $3$ & $4$ & $0$%
\end{tabular}%
\end{array}%
\end{array}%
\end{equation*}%
If $m=1,$ i.e. subjects may drop out in the final period only, then each
array with the last row deleted gives the corresponding minimal design $%
d_{\min }$. It can be verified that the information matrices of $d_{\min }$
will have rank 2, 1 and 4, respectively. So, as noted by Low, Lewis and
Prescott (1999), the minimal design corresponding to $d_{2plan}$ is
disconnected; indeed the only estimable direct treatment contrast in $%
d_{2\min }$ is $\tau _{0}-\tau _{1}+\tau _{2}-\tau _{3}.$ On the other hand
the minimal designs corresponding to $d_{1plan}$ and $d_{3plan}$ are
connected, the former has a nonzero eigenvalue $0.125$ with multiplicity $2;$
the nonzero eigenvalues of the latter are $2.61$ and $3.73,$ each with
multiplicity $2.$

\medskip

Since UBRMDs with three periods have been studied in Jones and Kenward
(2003) and Low (1995), henceforth we consider $t\geq 4$. The following lemma
shows that $d_{\min }$ corresponds to the maximal loss of information. For
nonnegative definite matrices $A$ and $B$ we use $A\succeq B$ to denote the
fact that $A-B$ is a nonnegative definite matrix, the L\"{o}wner order.

\medskip

\noindent \textbf{Lemma 2} $C_{D}^{d_{plan}}\succeq C_{D}^{d_{imp}}\succeq
C_{D}^{d_{\min }}.$

\medskip

This follows from known general results for linear models; for instance, it
is a consequence of Theorem 2.1 of Hedayat and Majumdar (1985). Lemma 2 says
that $d_{\min }$ has the "smallest information matrix" among all
possibilities for $d_{\mathrm{imp}}.$ The matrices, $C_{11}$, $C_{22}$ and $%
C_{12}$ for $d_{\min }$ are as given in Theorem 4. First, we need some
notation.

\noindent \textbf{The P and U matrices. }For $j=1,...,t,$ $P_{j}$ denotes a $%
t\times s$ matrix with $(h,i)$ entry $1$ if subject $i$ receives treatment $%
h $ in period $j$ of $d_{\mathrm{plan}}$; it is $0$ otherwise. For $%
j=0,1,...,t-1;$ $k=0,1,...,t-1,$ let, $U_{jk}=P_{t-j}P_{t-k}^{\prime }.$
Note that, since $d_{\mathrm{plan}}$ is a UBRMD,
\begin{equation}
P_{j}\mbox{\bf{1}}_{s}=g\mbox{\bf{1}}_{t},P_{j}^{\prime }\mbox{\bf{1}}_{t}=%
\mbox{\bf{1}}_{s}.
\end{equation}%
Also, the entries of $U_{jk}$ are nonnegative with row and column sums equal
to $g.$ Hence $\frac{1}{g}U_{jk}$ is a doubly stochastic matrix; in
particular $U_{jj}=gI_{t}.$ The following lemma gives the properties of
these matrices that we need.

\medskip

\noindent \textbf{Lemma 3} \ \textit{If }$U_{1},...,U_{M}$\textit{\ (not
necessarily distinct) are }$t\times t$\textit{\ matrices such that }$\frac{1%
}{g}U_{i}$\textit{\ is doubly stochastic for each }$i=1,...,M$\textit{, and }%
$a_{1},...,a_{M}$\textit{\ are nonnegative real numbers, then for }$x$%
\textit{\ }$\in $\textit{\ }$R^{t}$\textit{\ with }$x^{\prime }x=1,$\textit{%
\ }$x^{\prime }U_{i}x\leq g$\textit{\ for }$i=1,...,M,$\textit{\ and }$%
x^{\prime }(\sum a_{i}U_{i})^{\prime }(\sum a_{i}U_{i})x\leq g^{2}(\sum
a_{i})^{2}.$

\medskip

\noindent \textbf{Proof} \ If we write $W_{i}=\frac{1}{g}U_{i},$ then it
follows from the properties of doubly stochastic matrices (see, for example,
Bapat and Raghavan (1997, Chapter 2)) that $x^{\prime }W_{i}x\leq 1$. Also, $%
x^{\prime }(\sum a_{i}W_{i})^{\prime }(\sum a_{i}W_{i})x=x^{\prime }\left(
\sum \sum a_{i}a_{j}W_{i}^{\prime }W_{j}\right) x\leq \sum \sum a_{i}a_{j}%
\sqrt{x^{\prime }W_{i}^{\prime }W_{i}xx^{\prime }W_{j}^{\prime }W_{j}x}\leq
\sum \sum a_{i}a_{j}.$ The lemma follows.$\blacksquare $

\medskip

\noindent \textbf{Theorem 4} \textit{Let }$d_{\mathrm{plan}}$\textit{\ be a
UBRMD with }$t$\textit{\ treatments, }$t$\textit{\ time periods, }$s=gt$%
\textit{\ subjects, and let }$d_{\min }$\textit{\ consist of the first }$t-m$%
\textit{\ periods of }$d_{\mathrm{plan}}$\textit{. Then for }$d_{\min }$%
\textit{,}

\textit{(i) the information matrix for estimating direct and carry-over
treatment effects is given by~(3) and (4) with}%
\begin{equation}
C_{11}=\frac{g[(t-m)^{2}-m]}{t-m}I_{t}-\frac{g(t-2m)}{t-m}J_{t}-\frac{1}{t-m}%
\underset{j\neq k=0,...,m-1}{\sum \sum }U_{jk},
\end{equation}%
\begin{equation}
\begin{array}{c}
C_{22}=\frac{g}{t-m}[((t-m)^{2}-(t+1))I_{t}- \\
t^{-1}((t-m)^{2}-(t+1)-m(m+1))]J_{t}-\frac{1}{t-m}\underset{j\neq k=0,...,m}{%
\sum \sum }U_{jk},%
\end{array}%
\end{equation}%
\begin{equation}
C_{12}=\frac{g}{t-m}[(m+1)J_{t}-tI_{t}]-\underset{j=0}{\overset{m-1}{\sum }}%
U_{j(j+1)}-\frac{1}{t-m}\underset{j\neq k}{\overset{m-1}{\underset{j=0}{\sum
}}\underset{k=0}{\overset{m}{\sum }}}U_{jk}
\end{equation}

\textit{and (ii) if}%
\begin{equation}
t\geq 2m+2
\end{equation}%
\textit{then a g-inverse of }$C_{22}$\textit{\ is }$C_{22}^{-}=A^{-1},$%
\textit{\ where}%
\begin{equation}
A=\frac{g}{t-m}[(t-m)^{2}-(t+1)]I_{t}-\frac{1}{t-m}\underset{j\neq k=0,...,m}%
{\sum \sum }U_{jk}.
\end{equation}

\medskip

\noindent \textbf{Proof} \textit{(i)} Since $d_{\mathrm{plan}}$ is a \textit{%
UBRMD,} every treatment appears $s/t=g$ times in every period and $s$ times
in total. Also, in the order of application, each treatment is preceded by
every other treatment the same number of times and is not preceded by
itself. It follows that for the design $d_{\min },$%
\begin{eqnarray*}
N_{DS} &=&J_{t\times s}-\underset{j=0}{\overset{m-1}{\sum }}P_{t-j},\text{ }%
N_{CS}=J_{t\times s}-\underset{j=0}{\overset{m}{\sum }}P_{t-j}, \\
N_{DP} &=&gJ_{t\times (t-m)},\text{ }N_{CP}~=~\left[ \mathbf{0}%
_{t}~gJ_{t\times (t-m-1)}\right] , \\
r_{D} &=&g(t-m)\mbox{\bf{1}}_{t},\text{ }r_{C}=g(t-m-1)\mbox{\bf{1}}_{t}.
\end{eqnarray*}%
where $\mathbf{0}_{t}$ is a vector with $t$ zero elements. Inserting the
above formulae into~(\ref{eqGENinf}) and using (5) yields expressions $%
C_{11} $, $C_{22}$ and $C_{12}$ as in the statement of the theorem, after
some algebra.

\textit{(ii)} Using the relation $C_{22}\mbox{\bf{1}}_{t}=\mathbf{0}_{t},$
it can be verified that $C_{22}A^{-1}C_{22}=C_{22}$, as long as $A^{-1}$
exists; hence $A^{-1}$ is a g-inverse of $C_{22}$. We now show that
condition (9) guarantees the nonsingularity of $A.$ It follows from the fact
that the row sums of $\underset{j\neq k=0,...,m}{\sum \sum }U_{jk}$ are $%
gm(m+1),$ and Lemma 3, that the minimum eigenvalue of $A$ is%
\begin{equation}
\lambda _{\min }(A)=\frac{g}{t-m}[(t-m)^{2}-(t+1)-m(m+1)].
\end{equation}%
This is positive if $(t-m)^{2}-(t+1)-m(m+1)>0,$ which is equivalent to $%
t\geq 2m+2.$ $\blacksquare $

From the proof of Theorem 4, (9) guarantees that $A$ in (10) is positive
definite, and hence it is sufficient for the nonsingularity of $A$. This
condition plays a critical role in the derivation of bounds for the
eigenvalues of $C_{D}^{d_{\min }}.$ Also, it can be shown that (9) is a
necessary and sufficient condition for $rank(C_{22})=t-1.$ The next result
gives a bound on the eigenvalues of $C_{D}^{d_{\min }}$ which is used to
study the loss of precision for the estimators of the treatment contrasts
and the efficiency of $d_{\min }.$

\medskip

\noindent \textbf{Theorem 5} \textbf{\ }\textit{Suppose }$m\geq 1$\textit{\
and }$t\geq 2m+2$\textit{\ and }$d_{\mathrm{plan}}$\textit{\ is a UBRMD with
}$t$\textit{\ treatments, }$t$\textit{\ time periods, and }$s=gt$\textit{\
subjects. Denote the eigenvalues of }$C_{D}^{d_{\min }}$\textit{\ by }$%
g\theta _{0}=0,g\theta _{1},...,g\theta _{t-1}.$\textit{\ \ For }$%
r=1,...,t-1,$\textit{\ }$\theta _{r}\geq \theta _{L}(t,m),$\textit{\ where}%
\begin{equation}
\theta _{L}(t,m)=\frac{t}{t-m}\left[ (t-2m)-\frac{t(m+1)^{2}}{%
(t-m)^{2}-(t+1)-m(m+1)}\right] .
\end{equation}

\medskip

\noindent \textbf{Proof} \ \ Suppose $x\in R^{t},$ with $x^{\prime }x=1,$
and $x^{\prime }\mbox{\bf{1}}_{t}=0,$ such that $C_{D}^{d_{\min }}x=g\theta
_{r}x,$ $r=1,...,t-1.$ Then $g\theta _{r}=x^{\prime }C_{11}x-x^{\prime
}C_{12}A^{-1}C_{21}x.$ The maximum eigenvalue of $A^{-1}$ is $1/\left(
\lambda _{\min }(A)\right) $ where $\lambda _{\min }(A)$ is given by (11).
Hence, $g\theta _{r}\geq x^{\prime }C_{11}x-\frac{1}{\lambda _{\min }(A)}%
x^{\prime }C_{12}C_{21}x.$ If we write,%
\begin{equation*}
V=(t-m+1)\overset{m-1}{\underset{j=0}{\sum }}U_{j(j+1)}+\underset{k\neq j,j+1%
}{\overset{m-1}{\underset{j=0}{\sum }}\overset{m}{\underset{k=0}{\sum }}}%
U_{jk},
\end{equation*}%
then $C_{21}=-\frac{1}{t-m}\left[ gtI_{t}+V^{\prime }-g(m+1)J_{t}\right] .$
So,%
\begin{equation}
g\theta _{r}\geq x^{\prime }C_{11}x-\frac{1}{\lambda _{\min }\left( A\right) %
\left[ (t-m)^{2}\right] }\left[ (gt)^{2}+gtx^{\prime }(V+V^{\prime
})x+x^{\prime }VV^{\prime }x\right] .
\end{equation}%
The following inequalities can be derived by applying Lemma 3:%
\begin{eqnarray}
x^{\prime }C_{11}x &\geq &\frac{gt}{t-m}\left( t-2m\right)  \notag \\
x^{\prime }(V+V^{\prime })x &\leq &2g\left[ (t-m+1)m+m(m-1)\right] =2gmt
\notag \\
x^{\prime }VV^{\prime }x &\leq &\left( gmt\right) ^{2}.
\end{eqnarray}%
Inserting them into (13) and using the fact $\lambda _{\min }(A)>0,$ which
follows from the condition $t\geq 2m+2,$ we get a lower bound to $\theta
_{r} $ which, upon simplification, reduces to (12).$\blacksquare $

The results of Theorem 5 are useful in studying properties of $d_{\min }.$
Connectedness is a basic property of a design. A sufficient condition for $%
d_{\min }$ to be connected for direct treatment effects is $\theta
_{L}(t,m)>0.$ It follows from Lemma 2 that $d_{\mathrm{imp}}$ is connected
whenever $d_{\min }$ is connected. Corollary 6 follows from (12).

\medskip

\noindent \textbf{Corollary 6} \ \textit{Suppose }$d_{\mathrm{plan}}$\textit{%
\ is a UBRMD with }$t$\textit{\ treatments, }$t$\textit{\ time periods, and }%
$s=gt$\textit{\ subjects. A sufficient condition for the minimal design }$%
d_{\min }$\textit{\ to be connected is that }%
\begin{equation}
(t-2m)\left[ (t-m)^{2}-(t+1)-m(m+1)\right] -t(m+1)^{2}>0.
\end{equation}

\medskip

For a given $m,$ it follows from (15) that $d_{\min }$ is connected if a
polynomial in $t$ of degree $3$, which has leading coefficient $1$, is
positive. For $m=1,$ (15) reduces to $t^{3}-5t^{2}+4>0,$ i.e., $t\geq 5,$
and for $m=2$ it reduces to $t^{3}-9t^{2}+8t+12>0,$ i.e., $t\geq 8.$ These
observations lead to the following result.

\medskip

\noindent \textbf{Corollary 7} \ \textit{Suppose }$d_{\mathrm{plan}}$\textit{%
\ is a UBRMD with }$t$\textit{\ treatments, }$t$\textit{\ time periods, and }%
$s=gt$\textit{\ subjects. For each }$m\geq 1,$\textit{\ there is a positive
integer }$t^{\ast }(m)$\textit{\ such that the design }$d_{\min }$\textit{\
is connected if }$t\geq t^{\ast }(m).$\textit{\ In particular, for }$m=1,$%
\textit{\ }$d_{\min }$\textit{\ is connected whenever }$t\geq 5;$\textit{\
for }$m=2,$\textit{\ }$d_{\min }$\textit{\ is connected whenever }$t\geq 8.$

\medskip

One way to measure the \textit{goodness} of a connected design $d$ is by the
harmonic mean of the eigenvalues of the information matrix $C_{D}^{d},$ $%
H_{d}=(t-1)/trace(C_{D}^{d})^{+}$, where $C^{+}$ denotes the Moore-Penrose
inverse of $C.$ This is the value of the \textit{A-criterion}; hence $H_{d}$
is a measure of the precision of estimators of the direct treatment
contrasts for the design $d.$ It follows from Lemma 2 that $H_{d_{\min
}}\leq H_{d_{\mathrm{imp}}}\leq H_{d_{\mathrm{plan}}}.$ Since $d_{\mathrm{%
plan}}$ is the design that was chosen at the start of the experiment on the
basis of its desirable properties, especially efficiency, it is of interest
to examine the \textit{loss of precision} in the implemented design $d_{%
\mathrm{imp}}$\ with respect to $d_{\mathrm{plan}}$ due to subject dropout.
This loss may be measured by
\begin{equation*}
L_{d_{\mathrm{imp}}:d_{\mathrm{plan}}}=\frac{H_{d_{\mathrm{plan}}}-H_{d_{%
\mathrm{imp}}}}{H_{d_{\mathrm{plan}}}}=1-\frac{trace(C_{D}^{d_{\mathrm{plan}%
}})^{+}}{trace(C_{D}^{d_{\mathrm{imp}}})^{+}}.
\end{equation*}%
Clearly, the \textit{maximum loss of precision due to subject dropout} for $%
d_{\mathrm{plan}}$\ is given by%
\begin{equation*}
ML_{d_{\mathrm{plan}}}=L_{d_{\mathrm{\min }}:d_{\mathrm{plan}}}=\frac{H_{d_{%
\mathrm{plan}}}-H_{d_{\min }}}{H_{d_{\mathrm{plan}}}}=1-\frac{%
trace(C_{D}^{d_{\mathrm{plan}}})^{+}}{trace(C_{D}^{d_{\min }})^{+}},
\end{equation*}%
i.e., $ML_{d_{\mathrm{plan}}}\geq L_{d_{\mathrm{imp}}:d_{\mathrm{plan}}}.$

When $d_{\mathrm{plan}}$ is a UBRMD, we get using Theorem 5,%
\begin{equation*}
trace(C_{D}^{d_{\min }})^{+}=\underset{r=1}{\overset{t-1}{\sum }}\frac{1}{%
g\theta _{r}}\leq \frac{t-1}{g\theta _{L}(t,m)}.
\end{equation*}%
From Hedayat and Afsarinejad (1978), we obtain for the UBRMD $d_{\mathrm{plan%
}}$,%
\begin{equation*}
C_{D}^{d_{\mathrm{plan}}}=\frac{gt(t-2)(t+1)}{t^{2}-t-1}[I_{t}-\frac{1}{t}%
J_{t,t}],\text{ }trace(C_{D}^{d_{\mathrm{plan}}})^{+}=\frac{(t-1)(t^{2}-t-1)%
}{gt(t-2)(t+1)}.
\end{equation*}%
Therefore we obtain the following result.

\medskip

\noindent \textbf{Corollary 8} \ \textit{Suppose }$d_{\mathrm{plan}}$\textit{%
\ is a UBRMD with }$t$\textit{\ treatments, }$t$\textit{\ time periods, and }%
$s=gt$\textit{\ subjects. An upper bound to the maximum loss of precision
due to subject dropout is given by }$ML_{d_{\mathrm{plan}}}\leq UML(t,m)$%
\textit{\ where,}%
\begin{equation}
UML(t,m)=1-\frac{(t^{2}-t-1)\theta _{L}(t,m)}{t(t-2)(t+1)},
\end{equation}%
\textit{with }$\theta _{L}(t,m)$\textit{\ given by (12).}

\medskip

For $m=1$ and $t\geq 5,$ the values of $UML(t,1)$ for selected values of $t$
are given in Table 1, where the planned design $d_{\mathrm{plan}}$ is a
UBRMD. Similarly for $m=2$ and $t\geq 8,$ the values of $UML(t,2)$ for
selected values of $t$ are given in Table 2. As one would expect, for fixed $%
m,$ the bounds decrease with $t$ and become reasonably small when $t$ is
considerably larger than $t^{\ast }(m).$ In general, (16) is conservative.
Hence, the prospect of subject dropout may not be a big concern when $t$ is
much larger than $t^{\ast }(m)$.

If the UBRMD $d_{\mathrm{plan}}$ is chosen to have certain combinatorial
structures, the bound $UML(t,m)$ can be improved. One such structure is
considered next.

\medskip

\noindent \textbf{Type }$\mathcal{W}_{m}$\textbf{\ UBRMD}. Suppose the
subjects of the UBRMD $d_{\mathrm{plan}}$ can be partitioned into $g$ sets
of $t$ subjects each such that, within each group, every treatment appears
once in each of the periods $t-m,t-m+1,...,t$ for fixed $m\geq 1.$ Then for $%
j,k=0,...,m,$ $j\neq k$,
\begin{equation*}
U_{jk}=P_{t-j}P_{t-k}^{\prime }=\overset{g}{\underset{l=1}{\sum }}\Pi _{jkl},
\end{equation*}%
where each $\Pi _{jkl}$ is a permutation matrix of order $t$ and $P_{i}$ is
defined in Section 2. If\ for each $j,k=0,...,m,$ $j\neq k$ and $l=1,..,g$,
the eigenvalue $1$ of $\Pi _{jkl}$ has multiplicity one, then we say the
UBRMD $d_{\mathrm{plan}}$ is of \textit{type} $\mathcal{W}_{m}$.

If $m\geq 2$ an UBRMD of type $\mathcal{W}_{m}$ is also of type $\mathcal{W}%
_{m-1}.$ Examples of UBRMDs of \textit{type} $\mathcal{W}_{t-1}$ are UBRMD's
that are cyclically generated, for instance the Williams Latin Squares and
pairs of Williams Latin Squares given in Families 1 and 3 of Hedayat and
Afsarinejad (1978), and the class of sequentially counterbalanced squares
described by Isaac, Dean and Ostrom (2001).

It is known that the eigenvalues of a permutation matrix $\Pi $ of order $t$
are the roots of unity, $e^{i\frac{2\pi r}{t}}=Cos(2\pi r/t)+i\func{Si}%
n(2\pi r/t),$ $r=0,1,...,t-1,$ \textit{unless} the permutation can be
factored into the \textit{product of two or more disjoint cycles}, in which
case the multiplicity of $1$ as an eigenvalue of $\Pi $ is larger than one
(see, for example, Davis (1979)). If we set $\psi _{r}=Cos(\frac{2\pi r}{t})$
then, for a UBRMD of \textit{type} $\mathcal{W}_{m},$ the eigenvalues of \ $%
\Pi _{jkl}+\Pi _{jkl}^{\prime }$ are $2\psi _{r},$ $r=0,1,...,t-1.$ Since $%
U_{jk}+U_{kj}=\overset{g}{\underset{l=1}{\sum }}\left( \Pi _{jkl}+\Pi
_{jkl}^{\prime }\right) ,$ $\mbox{\bf{1}}_{t}^{\prime }x=0$ and $x^{\prime
}x=1$ implies $x^{\prime }(U_{jk}+U_{kj})x\leq 2g\psi _{1}.$ Hence the
inequalities (14) may be replaced by
\begin{equation*}
\begin{array}{c}
x^{\prime }C_{11}x\geq \frac{gt}{t-m}\left[ \left( t-2m\right) +\frac{m(m-1)%
}{t}\left( 1-\psi _{1}\right) \right] \\
x^{\prime }(V+V^{\prime })x\leq 2gmt\psi _{1} \\
x^{\prime }VV^{\prime }x\leq \left( gmt\right) ^{2}.%
\end{array}%
\end{equation*}%
If we insert these inequalities into (13) we obtain the following result.

\medskip

\noindent \textbf{Theorem 9} \ \textit{Suppose }$d_{\mathrm{plan}}$\textit{\
is a UBRMD of type }$W_{m}$\textit{\ for fixed }$m\geq 1$\textit{, }$t\geq
2m+2.$\textit{\ Denote the eigenvalues of }$C_{D}^{d}$\textit{\ by }$g\theta
_{0}=0,g\theta _{1},...,g\theta _{t-1}.$\textit{\ For }$r=1,...,t-1,$\textit{%
\ }$\theta _{r}\geq \theta _{L}^{\ast }(t,m),$\textit{\ where }%
\begin{equation*}
\theta _{L}^{\ast }(t,m)=\frac{t}{t-m}\left[ (t-2m)+\frac{m(m-1)}{t}\left(
1-\psi _{1}\right) -\frac{t\left( 1+2\psi _{1}m+m^{2}\right) }{%
(t-m)^{2}-(t+1)-m(m+1)}\right]
\end{equation*}%
\textit{with }$\psi _{1}=Cos\frac{2\pi }{t}.$

\medskip

Since $\psi _{1}<1,$ $\theta _{L}^{\ast }(t,m)>\theta _{L}(t,m).$ Hence
replacing $\theta _{L}(t,m)$ by $\theta _{L}^{\ast }(t,m)$ in (16) gives a
sharper upper bound to the \textit{maximum loss of precision due to subject
dropout} when $d_{\mathrm{plan}}$ is a UBRMD of \textit{type} $\mathcal{W}%
_{m}$, i.e., $ML_{d_{\mathrm{plan}}}\leq UML^{\ast }(t,m)<UML(t,m),$ where%
\begin{equation*}
UML^{\ast }(t,m)=1-\frac{(t^{2}-t-1)\theta _{L}^{\ast }(t,m)}{t(t-2)(t+1)}.
\end{equation*}

For $m=1$ and $t\geq 5,$ the values of the upper bound $UML^{\ast }(t,1)$ to
the maximum loss of precision due to subject dropout when $d_{\mathrm{plan}}$
is a UBRMD of \textit{type} $\mathcal{W}_{m}$ for selected values of $t$ are
given in Table 1. For $m=2$ and $t\geq 8$ the values of $UML^{\ast }(t,2)$
for selected values of $t$ are given in Table 2.

We now consider the efficiency of $d_{\min },$ the design that corresponds
to the \textit{worst case scenario}. Let $\emph{D}(t,gt,t-m)$ denote the
class of all connected crossover designs (not necessarily uniform or
balanced) in $t-m$ periods and $gt$ subjects, based on $t$ treatments. $%
d_{\min }$ belongs to this class. Since for an arbitrary design $d\in \emph{D%
}(t,gt,t-m),$ $\ trace(C_{D}^{d})^{+}\geq (t-1)^{2}/traceC_{D}^{d},$ a lower
bound for $trace(C_{D}^{d})^{+}$ may be obtained from an upper bound of $%
traceC_{D}^{d}.$ The latter bound can be obtained from Theorem 3 of Hedayat
and Yang (2004) (which generalized a bound of Stufken (1991)) as follows,
\begin{equation*}
\underset{d\in \emph{D}(t,gt,t-m)}{Max}traceC_{D}^{d}=gt(t-m-1)-\frac{%
2(gt-\delta ^{\ast })}{t-m}-\frac{(t-m-1)\delta ^{\ast 2}}{g(t-m)(t(t-m)-t-1)%
},
\end{equation*}%
where $\delta ^{\ast }$ is the nearest integer to $\frac{g(t(t-m)-t-1)}{t-m-1%
}.$ Since the choice $\delta ^{\ast }=$ $\frac{g(t(t-m)-t-1)}{t-m-1}$ gives
an upper bound to the maximum, it can be shown that $traceC_{D}^{d}\leq
gMTr(t,m),$ where%
\begin{equation*}
MTr(t,m)=t(t-m-1)-\frac{t(t-m-1)+1}{(t-m)(t-m-1)}.
\end{equation*}%
A measure of the efficiency of $d_{\min }$ in $\emph{D}(t,gt,t-m)$ is%
\begin{equation*}
EFF_{\emph{D}(t,gt,t-m)}^{d_{\min }}=\frac{\underset{d\in \emph{D}(t,gt,t-m)}%
{Min}trace(C_{D}^{d})^{+}}{trace(C_{D}^{d_{\min }})^{+}}\geq \frac{(t-1)^{2}%
}{gMTr(t,m)(trace(C_{D}^{d_{\min }})^{+})}.
\end{equation*}%
It follows from Theorems 5 and 9 that, if we define

\begin{equation*}
EL(t,m)=\frac{(t-1)\theta _{L}(t,m)}{MTr(t,m)}\text{ and }EL^{\ast }(t,m)=%
\frac{(t-1)\theta _{L}^{\ast }(t,m)}{MTr(t,m)},
\end{equation*}%
then the inequalities $EFF_{\emph{D}(t,gt,t-m)}^{d_{\min }}>EL(t,m)$ and $%
EFF_{\emph{D}(t,gt,t-m)}^{d_{\min }}>EL^{\ast }(t,m)$ give \textit{lower
bounds to the efficiency} of $d_{\min }$ when $d_{\mathrm{plan}}$ is a
general UBRMD and a UBRMD of \textit{type} $\mathcal{W}_{m},$ respectively.
Note that both $EL(t,m)$ and $EL^{\ast }(t,m)$ take values in $(0,1).$

For $m=1$ and $t\geq 5,$ the values of the lower bounds $EL(t,1)$ and $%
EL^{\ast }(t,1)$ to the efficiency of $d_{\min }$ in $\emph{D}(t,gt,t-1)$
for selected values of $t$ are given in Table 1. For $m=2$ and $t\geq 8,$
the values of $EL(t,2)$ and $EL^{\ast }(t,2)$ for selected values of $t$ are
given in Table 2. Since $EL(t,m)$ (or $EL^{\ast }(t,m)$) measures the
efficiency of $d_{\min }$ over \textit{all designs} in $\emph{D}(t,gt,t-m),$
not just those that are derived from UBRMDs, high values of this efficiency
bound suggest that a different starting design, instead of the UBRMD $d_{%
\mathrm{plan}},$ would not have resulted in a substantially better $d_{\min
}.$ An UBRMD $d_{\mathrm{plan}}$ that has a small value of $UML(t,m)$ (or $%
UML^{\ast }(t,m)$) and a large value of $EL(t,m)$ (or $EL^{\ast }(t,m)$) for
$d_{\min }$ clearly is a good design for use when there is a possibility of
subject dropout.

\section{Further results for the case $m=1$}

In the previous section we derived upper bounds $UML(t,m)$ and $UML^{\ast
}(t,m)$\ to the maximum loss of precision due to subject dropout $ML_{d_{%
\mathrm{plan}}}$. In this section, we first establish formulae for $ML_{d_{%
\mathrm{plan}}}$for two special families of UBRMDs. Then we indicate how to
select a starting design UBRMD $d_{\mathrm{plan}}$ for which $ML_{d_{\mathrm{%
plan}}}$ is small. For simplicity, we focus on the case $m=1,$ i.e.,
subjects remain in the study at least through period $t-1.$ We start with
definitions of the families of UBRMDs that we will study.

\medskip

\noindent \textbf{Class A type }$\mathcal{W}_{1}$\textbf{\ UBRMD: }Start
with a $t\times t$ square $W$\ that is a UBRMD with columns denoting
treatment sequences and rows denoting periods. Let $P_{t}(W)$ and $%
P_{t-1}(W) $ be the $t\times t$ matrices defined in Section 2 corresponding
to periods $t$ and $t-1,$ respectively, for square $W$, i.e., the $(h,i)$
entry of $P_{j}(W)$ is $1$ if the $(j,i)$ entry of $W$ is $h$; it is $0$
otherwise. Let $\Pi =P_{t}(W)P_{t-1}(W)^{\prime }.$ If $1$ is an eigenvalue
of $\Pi $ of multiplicity one, then the design $d_{\mathrm{plan}}$ that
assigns $g$ subjects to each sequence (column) of $W$ is called a Class A
type $\mathcal{W}_{1}$ UBRMD.

\medskip

\noindent \textbf{Class B type }$\mathcal{W}_{1}$\textbf{\ UBRMD: }We start
with two $t\times t$ squares $W_{1}$ and $W_{2}$ such that the $t\times 2t$
design $(W_{1}$ $W_{2})$ is a UBRMD with columns denoting treatment
sequences and rows denoting periods. For $\delta =1,2,$ let $P_{t}(W_{\delta
})$ and $P_{t-1}(W_{\delta })$ be the $t\times t$ matrices defined as in the
previous paragraph, and take $\Pi _{\delta }=P_{t}(W_{\delta
})P_{t-1}(W_{\delta })^{\prime }.$ Suppose $1$ is an eigenvalue of $\Pi
_{\delta }$ of multiplicity one for each $\delta =1,2.$ Suppose also that $%
W_{1}$ and $W_{2}$ are \textit{complementary} in the sense $\Pi _{2}=\Pi
_{1}^{\prime }.$ Then the design $d_{\mathrm{plan}}$ that assigns $g/2$
subjects to each sequence (column) of $W$ is called a Class B type $\mathcal{%
W}_{1}$ UBRMD.

Note that, for ease of implementation of the study, the experimenter will
generally assign several subjects to each of a small number of treatment
sequences (see Jones and Kenward (2003), p 159). Examples of Class A type $%
\mathcal{W}_{1}$ UBRMD when $t$ is even are the Williams Latin squares given
in Family 1 of Hedayat and Afsarinejad (1978) with $g$ subjects assigned to
each sequence, and examples of Class B type $\mathcal{W}_{1}$ UBRMD when $t$
is odd are the pair of William squares given in Family 3 of Hedayat and
Afsarinejad (1978) with $g/2$ subjects assigned to each sequence.

\medskip

\noindent \textbf{Theorem 10} \ \textit{Suppose }$t\geq 4.$\textit{\ (i) If }%
$d_{\mathrm{plan}}$\textit{\ is a Class A type }$W_{1}$\textit{\ UBRMD then,
for the minimal design }$d_{\min }$\textit{\ that consists of the first }$%
t-1 $\textit{\ periods of }$d_{\mathrm{plan}},$\textit{\ the eigenvalues of }%
$C_{D}^{d_{\min }}$\textit{\ are }$g\theta _{0}=0$\textit{\ and }$g\theta
_{r},$\textit{\ where}%
\begin{equation}
\theta _{r}=\frac{t}{t-1}\left[ t-2-\frac{2t(1+Cos(\frac{2\pi r}{t}))}{%
t(t-3)-2Cos(\frac{2\pi r}{t})}\right] ,\text{ }r=1,...,t-1.
\end{equation}

\textit{(ii) If }$d_{\mathrm{plan}}$\textit{\ is a Class B type }$W_{1}$%
\textit{\ UBRMD then, for the design }$d_{\min }$\textit{\ that consists of
the first }$t-1$\textit{\ periods of }$d_{\mathrm{plan}},$\textit{\ the
eigenvalues of }$C_{D}^{d_{\min }}$\textit{\ are }$g\theta _{0}=0$\textit{\
and }$g\theta _{r},$\textit{\ where}%
\begin{equation}
\theta _{r}=\frac{t}{t-1}\left[ t-2-\frac{t(1+Cos(\frac{2\pi r}{t}))^{2}}{%
t(t-3)-2Cos(\frac{2\pi r}{t})}\right] ,\text{ }r=1,...,t-1.
\end{equation}

\medskip

\noindent \textbf{Proof} \ Write $U=U_{01}=P_{t}P_{t-1}^{\prime }.$ Then $%
C_{11}=(gt(t-2)/(t-1))(I_{t}-\frac{1}{t}J_{t}),$ \ $%
C_{12}=-(1/(t-1))(gtI_{t}-2gJ_{t}+tU),$ and $%
C_{22}=A-(g(t^{2}-3t-2)/(t(t-1)))J_{t},$ where $A=((gt(t-3)/(t-1))I_{t}-%
\frac{1}{t-1}(U+U^{\prime }).$ Consider the spectral decomposition, $%
U+U^{\prime }=\underset{r=0}{\overset{t-1}{\sum }}\alpha
_{r}h_{r}h_{r}^{\prime },$ with $h_{r}^{\prime }h_{r}=1,$ $r=0,1,...,t-1,$ $%
h_{r}^{\prime }h_{q}=0,$ for $r\neq q;$ $\alpha _{0}=2g,$ $h_{0}=\frac{1}{%
\sqrt{t}}\mbox{\bf{1}}_{t}.$ For $r=1,...,t-1$, $h_{r}^{\prime }\mbox{\bf{1}}%
_{t}=0.$ Let $\gamma _{r}=(I_{t}+\frac{1}{g}U)h_{r},$ $r=0,1,...,t-1.$ It
can be shown that%
\begin{equation*}
C_{D}^{d_{\min }}=\frac{gt}{t-1}[(t-2)(I_{t}-\frac{1}{t}J_{t,t})-gt\overset{%
t-1}{\underset{r=1}{\sum }}(gt(t-3)-\alpha _{r})^{-1}\gamma _{r}\gamma
_{r}^{\prime }].
\end{equation*}%
Note that, by Lemma 3, $\alpha _{r}\leq 2g.$ Hence for $t\geq 4,$ $%
gt(t-3)-\alpha _{r}\geq gt(t-3)-2g>0.$

\textit{(i)} In this case, $U=g\Pi .$ For $r=1,...,t-1,$ $(U+U^{\prime
})h_{r}=\alpha _{r}h_{r}$ implies $\alpha _{r}=2g\psi _{r},$where $\psi
_{r}=Cos(\frac{2\pi r}{t}).$ Since $\gamma _{r}=\left( I_{t}+\Pi \right)
h_{r},$ $\gamma _{l}^{\prime }\gamma _{r}=h_{l}^{\prime }(2I_{t}+\Pi +\Pi
^{\prime })h_{r}=\left( 2+2\psi _{r}\right) h_{l}^{\prime }h_{r}.$
Therefore, if we write $\gamma _{0}^{\ast }=\frac{1}{\sqrt{t}}\mbox{\bf{1}}%
_{t},$ $\gamma _{r}^{\ast }=\left( 2+2\psi _{r}\right) ^{-1/2}\gamma
_{r}=\left( 2+2\psi _{r}\right) ^{-1/2}\left( I_{t}+\Pi \right) h_{r},$ $%
r=1,...,t-1,$ then $\{\gamma _{0}^{\ast },\gamma _{1}^{\ast },...,\gamma
_{t-1}^{\ast }\}$ is an orthogonal and normalized basis of $\mathcal{R}^{t}$%
, and for $r=1,...,t-1,$
\begin{equation*}
C_{D}^{d_{\min }}\gamma _{r}^{\ast }=\frac{gt}{t-1}\left[ t-2-\frac{t\left(
2+2\psi _{r}\right) }{t(t-3)-2\psi _{r}}\right] \gamma _{r}^{\ast }=\frac{gt%
}{t-1}\left[ t-2-\frac{2t(1+Cos(\frac{2\pi r}{t}))}{t(t-3)-2Cos(\frac{2\pi r%
}{t})}\right] \gamma _{r}^{\ast },
\end{equation*}%
which establishes \textit{(i)}.

\textit{(ii)} Here $U=\frac{g}{2}(\Pi _{1}+\Pi _{1}^{\prime })=U^{\prime }.$%
\ For $r=1,...,t-1,$ $(U+U^{\prime })h_{r}=\alpha _{r}h_{r}$ implies $\alpha
_{r}=2g\psi _{r};$ $\gamma _{r}=\left( I_{t}+\frac{1}{g}U\right)
h_{r}=(1+\psi _{r})h_{r}.$ It follows that, $h_{0},h_{1},...,h_{r}$ are
orthogonal and normalized eigenvectors of $C_{D}^{d_{\min }}$, and%
\begin{equation*}
C_{D}^{d_{\min }}h_{r}=\frac{gt}{t-1}\left[ t-2-\frac{t(1+Cos(\frac{2\pi r}{t%
}))^{2}}{t(t-3)-2Cos(\frac{2\pi r}{t})}\right] h_{r},
\end{equation*}%
for $r=1,...,t-1.$ This establishes \textit{(ii)}. $\blacksquare $

For $t=4$ it can be shown that $\theta _{1}=0,$ $\theta _{2}=2.67,$ and $%
\theta _{3}=0.$ Hence $d_{\min }$ is not connected. For $t\geq 5$ however,
it follows from Corollary 7 that $d_{\min }$ is connected. The following
Corollary is immediate.

\medskip

\noindent \textbf{Corollary 11} \ \textit{Suppose }$t\geq 5.$\textit{\ If }$%
d_{\mathrm{plan}}$\textit{\ is a Class A or Class B type }$W_{1}$\textit{\
UBRMD then the loss of precision due to subject dropout is }%
\begin{equation*}
L_{d_{\mathrm{imp}}:d_{\mathrm{plan}}}\leq L_{d_{\mathrm{\min }}:d_{\mathrm{%
plan}}}=ML_{d_{\mathrm{plan}}}=1-\frac{(t-1)(t^{2}-t-1)}{t(t-2)(t+1)}\left(
\overset{t-1}{\underset{r=1}{\sum }}\frac{1}{\theta _{r}}\right) ^{-1},
\end{equation*}%
\textit{where the }$\theta _{r}$\textit{'s are given by (17) and (18) for
Class A and B, respectively.}

\medskip

Theorem 10 also gives a lower bound to the efficiency of $d_{\min }$ in $%
\emph{D}(t,gt,t-1),$ $EFF_{\emph{D}(t,gt,t-1)}^{d_{\min }}>EL^{AB}(t)$, where%
\begin{equation*}
EL^{AB}(t)=\frac{(t-1)^{2}}{MTr(t,1)}\left( \overset{t-1}{\underset{r=1}{%
\sum }}\frac{1}{\theta _{r}}\right) ^{-1}
\end{equation*}%
and where $MTr(t,1)=t(t-2)-\frac{t-1}{t-2}.$ The difference between the
bounds $EL^{AB}(t)$ and $EL^{\ast }(t,1)$ defined in Section 2 is that in
the former (which applies to Class A or Class B type $\mathcal{W}_{1}$
UBRMD) exact values of the eigenvalues $\theta _{r}$ are used while in the
latter (which applies to any type $\mathcal{W}_{1}$ UBRMD) these are
replaced by the lower bound $\theta _{L}^{\ast }(t,1)$.

Table 3 gives $ML_{d_{\mathrm{plan}}}$ and $EL^{AB}(t)$ for Class A and
Class B type $\mathcal{W}_{1}$ UBRMD for selected values of $t$. Note that
the values of the maximum loss due to subject dropout $ML_{d_{\mathrm{plan}%
}} $ are substantially lower than the upper bounds given in Table 1, while
the values of the efficiency bound $EL^{AB}(t)$ are higher.

For $g>1,$ a Class A type $\mathcal{W}_{1}$ UBRMD consists of $g$
replications of a square while a Class B type $\mathcal{W}_{1}$ UBRMD
consists of $g/2$ replications of a pair of squares. The next result implies
that the loss due to subject dropout may be smaller if distinct squares (or
pair of squares) are used instead of replications.

\medskip

\noindent \textbf{Corollary 12} \ \textit{Suppose }$t\geq 5.$\textit{\ If }$%
d_{\mathrm{plan}}$\textit{\ is the union of }$g$\textit{\ }$t\times t$%
\textit{\ Class A type }$W_{1}$\textit{\ UBRMDs \ or }$g/2$\textit{\ }$%
t\times 2t$\textit{\ Class B type }$W_{1}$\textit{\ UBRMDs that are not
necessarily identical, then }%
\begin{equation*}
L_{d_{\mathrm{imp}}:d_{\mathrm{plan}}}\leq L_{d_{\mathrm{\min }}:d_{\mathrm{%
plan}}}=ML_{d_{\mathrm{plan}}}\leq 1-\frac{(t-1)(t^{2}-t-1)}{t(t-2)(t+1)}%
\left( \overset{t-1}{\underset{r=1}{\sum }}\frac{1}{\theta _{r}}\right)
^{-1},
\end{equation*}%
\textit{where the }$\theta _{r}$\textit{'s are given by (17) and (18) for
Class A and B, respectively.}

\medskip

\noindent \textbf{Proof} \ We sketch a proof for Class A type $\mathcal{W}%
_{1}$ UBRMDs; the proof for Class B type $\mathcal{W}_{1}$ UBRMDs is
identical. Suppose $d_{\mathrm{\min }}=\underset{i=1}{\overset{g}{\dbigcup }}%
d_{i},$ where $d_{i}$ is the minimal design for a Class A type $\mathcal{W}%
_{1}$ UBRMD design based on $t$ subjects for $i=1,...,g$. It follows from
Theorem 2.1 of Hedayat and Majumdar (1985) that, $C_{D}^{d_{\mathrm{\min }%
}}\succeq \overset{g}{\underset{i=1}{\sum }}C_{D}^{d_{i}}.$ Take $%
B=C_{D}^{d}+\frac{1}{t}J_{t},$ and $B_{i}=C_{D}^{d_{i}}+\frac{1}{gt}J_{t},$
for $i=1,...,g.$ Since each $d_{i}$ is connected, the matrices $%
B_{1},...,B_{g}$ and $B$ are positive definite. Clearly, $B\succeq \overset{g%
}{\underset{i=1}{\sum }}B_{i}.$ It follows that
\begin{equation*}
B^{-1}\preceq \left( \overset{g}{\underset{i=1}{\sum }}B_{i}\right)
^{-1}\preceq \frac{1}{g^{2}}\left( \overset{g}{\underset{i=1}{\sum }}%
B_{i}^{-1}\right) ,
\end{equation*}%
where the first inequality is well known in matrix theory and the second is
given in Bapat and Raghavan ((1997), Theorem 3.11.1). \ It can be shown
that, $B^{-1}=\left( C_{D}^{d_{\mathrm{\min }}}\right) ^{+}+\frac{1}{t}%
J_{t}, $ and $B_{i}^{-1}=\left( C_{D}^{d_{i}}\right) ^{+}+\frac{g}{t}J_{t},$
for $i=1,...,g.$ Hence,$\left( C_{D}^{d_{\mathrm{\min }}}\right) ^{+}\preceq
\frac{1}{g^{2}}\overset{g}{\underset{i=1}{\sum }}\left( C_{D}^{d_{i}}\right)
^{+}.$ This implies,
\begin{equation*}
trace\left( C_{D}^{d_{\mathrm{\min }}}\right) ^{+}\leq \frac{1}{g^{2}}%
\overset{g}{\underset{i=1}{\sum }}trace\left( C_{D}^{d_{i}}\right) ^{+}=%
\frac{1}{g^{2}}\overset{g}{\underset{i=1}{\sum }}\underset{r=1}{\overset{t-1}%
{\sum }}\frac{1}{\theta _{r}}=\frac{1}{g}\underset{r=1}{\overset{t-1}{\sum }}%
\frac{1}{\theta _{r}}.
\end{equation*}%
The result follows.$\blacksquare $

For the setup of Corollary 12, it is clear that a lower bound to the
efficiency of $d_{\min }$ in $\emph{D}(t,gt,t-1)$ is
\begin{equation*}
EFF_{\emph{D}(t,gt,t-1)}^{d_{\min }}>(t-1)^{2}/\left( MTr(t,1)(\underset{r=1}%
{\overset{t-1}{\Sigma }}1/\theta _{r})\right) .
\end{equation*}

Corollary 12 indicates that the use of distinct Class A or Class B type $%
\mathcal{W}_{1}$ UBRMDs instead of replications of the same design will not
increase the maximum loss of precision due to subject dropout $ML_{d_{%
\mathrm{plan}}}.$ There are examples where $ML_{d_{\mathrm{plan}}}$actually
decreases. In their Example 2, Low, Lewis and Prescott (1999) studied the
case $t=4,$ $s=24,$ $m=1$ and showed that the use of distinct William
Squares instead of replications of the same square reduces $ML_{d_{\mathrm{%
plan}}}$. An example for $t=6$ is given below. The implication is that a
UBRMD with more distinct sequences is likely to perform better under subject
dropout.

Consider the "extreme" design $d_{\mathrm{plan}}^{e}$ that consists of one
subject assigned to each of the $t!$ possible sequences $(s=t!)$. For the
case $m=1,$ it can be shown that $U_{01}=((t-2)!)[J_{t}-I_{t}]$ and the
information matrix of the minimal design $d_{\min }^{e}$ is%
\begin{equation*}
C_{D}^{d^{e}}=\frac{at(t-2)[(t-2)!]}{t-1}\left( I_{t}-\frac{1}{t}%
J_{t}\right) ,\text{ where }a=\frac{t^{4}-5t^{3}+6t^{2}+t-2}{%
t^{3}-4t^{2}+3t+2}.
\end{equation*}%
For $d_{\mathrm{plan}}^{e}$ the maximum loss of precision due to subject
dropout is,
\begin{equation*}
ML_{d_{\mathrm{plan}}^{e}}=1-a(t^{2}-t-1)/((t-1)^{2}(t+1)).
\end{equation*}%
Numerical studies indicate that this is the smallest value of $ML_{d_{%
\mathrm{plan}}}$ among all UBRMDs with $t!$ or fewer subjects. We are
currently investigating the nature of planned designs that attain the
minimum and maximum values of $ML_{d_{\mathrm{plan}}}$, as well as designs
that fall in between these extremes. However, as mentioned earlier in this
section, a small number of treatment sequences is generally preferred, so it
is doubtful that crossover designs with a large number of distinct sequences
will be used widely in practice.

\medskip

\noindent \textbf{Example 13} \ Let $t=6$ and $s=12g_{0}$. The array below
consists of two Williams Latin squares.
\begin{equation*}
\begin{array}{cc}
\text{Square 1} & \text{Square 2} \\
\begin{array}{cccccc}
1 & 2 & 3 & 4 & 5 & 0 \\
0 & 1 & 2 & 3 & 4 & 5 \\
2 & 3 & 4 & 5 & 0 & 1 \\
5 & 0 & 1 & 2 & 3 & 4 \\
3 & 4 & 5 & 0 & 1 & 2 \\
4 & 5 & 0 & 1 & 2 & 3%
\end{array}
&
\begin{array}{cccccc}
2 & 5 & 1 & 3 & 0 & 4 \\
4 & 2 & 5 & 1 & 3 & 0 \\
5 & 1 & 3 & 0 & 4 & 2 \\
0 & 4 & 2 & 5 & 1 & 3 \\
1 & 3 & 0 & 4 & 2 & 5 \\
3 & 0 & 4 & 2 & 5 & 1%
\end{array}%
\end{array}%
\end{equation*}%
Suppose $d_{\mathrm{plan}}^{1}$ is a design that assigns $2g_{0}$ subjects
to the first six columns of the array and $d_{\mathrm{plan}}^{2}$ a design
that assigns $g_{0}$ subjects to each of the twelve columns. If $m=1,$ then
the maximum loss of precision due to subject dropout are $ML_{d_{\mathrm{plan%
}}^{1}}=0.30$ and $ML_{d_{\mathrm{plan}}^{2}}=0.24.$ The design $d_{\mathrm{%
plan}}^{e}$ can be constructed when $g_{0}=60.$ For this design, $ML_{d_{%
\mathrm{plan}}^{e}}=0.20.$

\medskip

Since estimation of the residual effects of the treatments is sometimes at
least a secondary focus of experiments, we conclude this section with a
brief consideration of the information matrix for the residual treatment
effects of the minimal design $d_{\min }$ obtained from a UBRMD $d_{\mathrm{%
plan}}$ when $m=1.$ Note that, $C_{R}^{d_{imp}}\succeq C_{R}^{d_{\min }}$.
Also, $C_{R}^{d_{\mathrm{plan}}}=C_{22}-C_{21}C_{11}^{-}C_{12}.$ Suppose $%
t\geq 3.$ \ Then, it can be shown that $\left( (t-1)/(gt(t-2))\right) I_{t}$
is a generalized inverse of $C_{11}.$ From (2) and (4) we obtain,%
\begin{eqnarray*}
C_{R}^{d_{\min }} &=&\left( \frac{gt(t^{2}-5t+5)}{(t-1)(t-2)}\right) \left[
I_{t}-\frac{1}{t}J_{t,t}\right] -\left( \frac{2}{t-2}\right) \left[
U+U^{\prime }\right] \\
&&-\left( \frac{t}{g(t-1)(t-2)}\right) U^{\prime }U+\left( \frac{g(5t-4)}{%
t(t-1)\left( t-2\right) }\right) J_{t,t}.
\end{eqnarray*}%
Using this, it can be shown that if the design $d_{\mathrm{plan}}$ is a
Class A or Class B type $\mathcal{W}_{1}$ UBRMD then $d_{\min }$ is
connected for residual treatment effects whenever $t=3,$ or $t\geq 5.$ For $%
t=4,$ Low, Lewis and Prescott (1999) have shown that if $d_{\mathrm{plan}}$
is a Williams Latin square then $d_{\min }$ is disconnected for the residual
treatment effects.

\medskip

\textbf{ACKNOWLEDGEMENT} The authors would like to thank two referees, an
associate editor and the Editor for insightful comments that improved the
manuscript substantially. Dibyen Majumdar's research is sponsored by NSF
grant DMS-0204532 and Angela Dean's research is sponsored by NSF grant
SES-0437251.

\begin{center}
\bigskip {\large References}
\end{center}

Bapat, R. B. and Raghavan, T. E. S. (1997). \textit{Nonnegative Matrices and
Applications}. Cambridge University Press.

Cheng, C.-S. and Wu, J. (1980). Balanced repeated measurement designs.
\textit{Ann. Statist}. \textbf{8}, 1272-1283. (Corrigendum. \textit{Ann.
Statist}. (1983), \textbf{11}, 349).

Davis, P. J. (1979). \textit{Circulant Matrices}, John Wiley and Sons, New
York.

Diggle, P. D. and Kenward, M. G. (1994). Informative dropout in longitudinal
data analysis (with discussion). \textit{J. Roy. Statist. Soc}. C \textbf{43}%
, 49-93.

Ghosh, S. (1979). On robustness of designs against incomplete data. \textit{%
Sankhya} \textbf{40}, 204-208.

Ghosh, S. (1982). Robustness of BIBD against the unavailability of data.
\textit{J. Statist. Plann. Inf.} \textbf{6}, 29-32.

Godolphin, J. D. (2004). Simple pilot procedures for the avoidance of
disconnected experimental designs. \textit{Appl. Statist}. \textbf{53},
133-147.

Higham, J. (1998). Row-complete latin squares of every composite order
exist. \textit{J. Combin Design} \textbf{6}, 63-77.

Hedayat, A. S. and Afsarinejad, K. (1978). Repeated measurements designs,
II. \textit{Ann. Statist}. \textbf{6}, 619-628.

Hedayat, A. S. and Majumdar, D. (1985). Combining experiments under
Gauss-Markov models. \textit{J. Amer. Statist. Assoc}. \textbf{80}, 698-703.

Hedayat, A. S. and Yang, M. (2003). Universal optimality of balanced uniform
crossover designs. \textit{Ann. Statist}. \textbf{31}, 978-983.

Hedayat, A. S. and Yang, M. (2004). Universal optimality of selected
crossover designs. \textit{J. American Statistical Association} \textbf{99},
461-466.

Isaac, P. D., Dean, A. M. and Ostrom, T. (2001). Generating pairwise
balanced Latin Squares. \textit{Statistics and Applications} \textbf{3},
25-46.

Jones B. and Kenward M. G. (2003). \textit{Design and Analysis of Cross-over
Trials}. CRC Press, London.

Kunert, J. (1984). Optimality of balanced repeated measurements designs.
\textit{Ann. Statist.} \textbf{12}, 1006-1017.

Low, J. L. (1995). The design of crossover designs subject to dropout. PhD
thesis, School of Mathematics, University of Southampton.

Low, J. L., Lewis, S. M. and Prescott, P. (1999). Assessing robustness of
crossover designs to subjects dropping out. \textit{Statistics and Computing}
\textbf{9}, 219-227.

Ratkowsky, D. A., Evans, M. A. and Alldredge, J. R. (1992). \textit{%
Crossover Experiments: Design, Analysis and Application}. Marcel Dekker, New
York.

Senn, S. (2002). \textit{Cross-over Trials in Clinical Research}. Wiley, New
York.

Stufken, J. (1991). Some families of optimal and efficient repeated
measurements designs, \textit{J. Statist. Plann. Inf.} \textbf{27}, 75-83.

Stufken, J. (1996). Optimal crossover designs. \textit{Design and Analysis
of Experiments. Handbook of Statistics} \textbf{13} (S. Ghosh and C. R. Rao,
eds.), 63-90, North Holland, Amsterdam.

Williams, E. J. (1949). Experimental designs balanced for the estimation of
residual effects of treatments. \textit{Austral. J. Sci. Res.} A2,
149-168.\bigskip

\pagebreak

\begin{center}
${\large TABLES}$

\bigskip

Table 1: Upper bounds to the maximum loss of precision due to subject dropout

and lower bound to the efficiency of the minimal design when $m=1$

\bigskip

\bigskip
\begin{tabular}{|c|c|c|c|c|c|c|}
\hline
$t$ & 5 & 6 & 7 & 8 & 9 & 10 \\ \hline
$UML(t,1)$ & 0.87 & 0.48 & 0.33 & 0.25 & 0.20 & 0.17 \\ \hline
$UML^{\ast }(t,1)$ & 0.64 & 0.40 & 0.30 & 0.23 & 0.19 & 0.16 \\ \hline
$EL(t,1)$ & 0.18 & 0.66 & 0.81 & 0.88 & 0.92 & 0.93 \\ \hline
$EL^{\ast }(t,1)$ & 0.49 & 0.76 & 0.85 & 0.90 & 0.93 & 0.94 \\ \hline
\end{tabular}

\bigskip

Table 2: Upper bounds to the maximum loss of precision due to subject dropout

and lower bound to the efficiency of the minimal design when $m=2$

\bigskip

\bigskip
\begin{tabular}{|c|c|c|c|c|c|c|}
\hline
$t$ & 8 & 9 & 10 & 11 & 12 & 16 \\ \hline
$UML(t,2)$ & 0.90 & 0.63 & 0.48 & 0.39 & 0.33 & 0.20 \\ \hline
$UML^{\ast }(t,2)$ & 0.81 & 0.59 & 0.46 & 0.38 & 0.32 & 0.21 \\ \hline
$EL(t,2)$ & 0.15 & 0.50 & 0.67 & 0.77 & 0.82 & 0.93 \\ \hline
$EL^{\ast }(t,2)$ & 0.27 & 0.55 & 0.69 & 0.78 & 0.83 & 0.94 \\ \hline
\end{tabular}

\bigskip

Table 3: Maximum loss of precision due to subject dropout

and lower bound to the efficiency of the minimal design

for Class A and Class B type $\mathcal{W}_{1}$ UBRMD $d_{\mathrm{plan}}$
when $m=1$

\bigskip

\bigskip
\begin{tabular}{|c|c|c|c|c|c|c|}
\hline
t & 5 & 6 & 7 & 8 & 9 & 10 \\ \hline
Class & B & A & B & A & B & A \\ \hline
$ML_{d_{\mathrm{plan}}}$ & 0.35 & 0.30 & 0.20 & 0.18 & 0.14 & 0.13 \\ \hline
$EL^{AB}(t)$ & 0.90 & 0.89 & 0.97 & 0.97 & 0.98 & 0.98 \\ \hline
\end{tabular}
\end{center}

\bigskip

\end{document}